# Four Cubes
Szymon Łukaszyk

A short survey on the properties of four graphs constructed in $\{0, 1\}^n$ Boolean space is presented. Flexible activation function of an artificial neuron in a sparse distributed memory model is defined on the basis of the Ugly duckling theorem. Cotan Laplacian on 2-face triangulation of $n$-cube has degenerate spectrum of eigenvalues corresponding to the Hamming distance distribution of $\{0, 1\}^n$ space. Degenerate spectrum of eigenvalues of the cotan Laplacian defined on the graph comprising $2^n$ 2-face triangulated $n$-cubes sharing common origin includes all integers from 0 to $3n$ without the eigenvalue of $3n$-1 (multiplicities of the same eigenvalues form A038717 OEIS sequence), while the multiplicities of the same eigenvalues $[-n\sqrt{2}, n\sqrt{2}]$ of the adjacency matrix of $2^n$-cube form trinomial triangle. The distance matrix of this graph, providing further OEIS sequences, as well as its relation with Buckminster Fuller vector equilibrium is also discussed.

## I. Introduction

All that we perceive is information, the resolution of uncertainty, or knowledge that signifies our understanding of a concept. The information is measured in bits, discrete basic units of information. On the other hand we perceive the universe in 3 spatial dimensions and unidirectional temporal one. But the arrow of time relates not only to philosophy and the 2nd law of thermodynamics but also to biological evolution.

However, the question is why do we sense such a dimensionality? Does it result from some kind of a minimum energy condition? Or could there exist individuals (flatlanders?) perceiving our (or their) universe in $n$ dimensions, with finite $n \ne 3$ *spatial* ones? Pseudo-Riemannian manifolds used in general relativity theory are not bound to four dimensions and the combinatorial proof of the Boltzmann's $H$-theorem [LB77] introducing the concept of energy quantization which led to the development of quantum theory [TK78], is irrelevant to any particular dimensionality of space, in which Ludwig Boltzmann considered the molecules to collide.

An attempt to answer this question and link the perceived dimensionality with space-independent information through perception of a biological entity and in an information theoretic approach was the subject of this paper. Although this goal has not been completed, a survey on four cubes is presented in a hope that any of those who, like the author, are lost in math and wander astray in search of physical hard to vary models of reality will find here anything handy for their research.

## II. The Cubes

A recurrence relation

$$f_n \doteq 2f_{n-2}/n$$

for $n \in \mathbb{N}_0$, where $f_0 := 1$ and $f_1 := 2$ allows to express the volumes and surfaces of all $n$-balls as

$$V_n = \pi^{\lfloor n/2 \rfloor} f_n R^n \quad \text{(VNB)}$$

$$S_n = n\pi^{\lfloor n/2 \rfloor} f_n R^{n-1} \quad \text{(SNB)}$$

$S_n = 2\pi R V_{n-2}$ and the subsequent exponents of $\pi^{\lfloor n/2 \rfloor}$ are even if $n = 4k$ or $n = 4k+1$ for $k \in \mathbb{N}_0$ and odd otherwise. The third unit cube of those that we shall now introduce is Eulerian circuitable also if $n = 4k-1$ or $n = 4k$ for $k \in \mathbb{N}$. The first two cubes are commonly known. The third one is probably also known but it has some, as author believes unknown, and interesting properties. The last one contains $2^n$ first, second or third ones.

### 1. *n*-cube

Every simplicial $n$-manifold inherits a natural topology from $\mathbb{R}^n$, while local and global invariants are easily discovered when expressed in discrete forms rather than by staring at the indices [DKT08] of some field equations. The key aspect of this this approach is that it disentangles the topological (metric-independent) and geometrical (metric-dependent) content of the modelled quantities, keeping their intrinsic structure intact. Operators, that in the continuous theory do not use metric information maintain this property in the discrete theory as well [AH03]. Thus the simplex formulation is equivalent to the continuous calculus. At least in this regard.

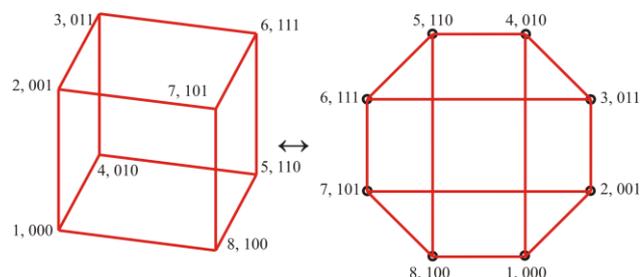

Fig. 1: Gray encoded 3-cube with its associated graph.

$n$-cube is the smallest proper $n$-dimensional hole that can be created in in a simplicial $n$-manifold. Removing just one $n$-simplex could not serve as a model of a continuous hole: Stokes integration over an $n$-manifold with one missing $n$-simplex is the same as if this $n$-simplex was present. Also an $n$-orthoplex (dual to $n$-cube) could not be a proper hole (for $n \geq 3$), since having $2^{k+1}(n:k+1)$ $k$-faces it formally does not have an $n$-face[1]. Besides, all the remaining faces of $n$-orthoplex (just like all the faces of $n$- simplex) are ($n-k$)-simplices, while all the faces of $n$-cube are ($n-k$)-cubes; $n$-orthoplex cannot be extruded from a point like an $n$-cube[2]. $n$-cube defines a Cartesian coordinate system of $\mathbb{R}^n$ or at least a $2^{-n}$ part of it (cf. Fig. 6).

---
[1] $(n:n+1) = 0$ if binomial coefficient is defined in terms of a falling factorial $(n)_{k+1}$.

If the $n$-tuples of the addresses $a_m$ of the vertices $m = 1, 2,\ldots, 2^n$ of a unit $n$-cube are ordered using Gray code, so that the subsequent addresses differ by just one bit, the bipartite graph associated to $n$-cube has a form shown in Fig. 1.

One can define a distance matrix $D$ for $n$-cube

$$D_{lm} \doteq d_{HM}(a_l, a_m) \qquad \text{(DMT)}$$

where $d_{HM}(a_l, a_m)$ is the Hamming distance between the addresses of the vertices $l$ and $m$.

Matrix $D$ is (the list is not exhaustive):
- symmetric and hollow (zeros on the main diagonal), so its trace is zero;
- $D_{lm} \leq D_{lk} + D_{km}$ for all $k$ (the triangle inequality);
- even dimensional ($2^n \times 2^n$) and its dimension is divisible by 4 for n $\geq$ 2;[3]
- if the subsequent vertices are binary encoded, entries with $d_{HM} = n$ reside on the counterdiagonal;
- if the subsequent vertices are Gray encoded, entries with $d_{HM} = 1$ reside on the counterdiagonal;
- if the subsequent vertices are Gray or binary encoded $D$ is also centrosymmetric.

Matrix $D$ leads us to {$n$}-cube having edges carrying Hamming distances between particular vertices.

## 2. {$n$}-cube

{$n$}-cube shown in Fig. 2 is just a fancy name for Boolean $\{0, 1\}^n$ address space. It is easy to see that {$n$}-cube is isomorphic to ($2^n-1$)-simplex. On the other hand a regular $n$-simplex can be inscribed in $n$-cube if and only if $n = 2^m - 1$ for $m \in \mathbb{N}_0$[4]. The degree of any vertex $m$ of {$n$}-cube is $2^{1-n}(2^n : 2)$, which is odd number.

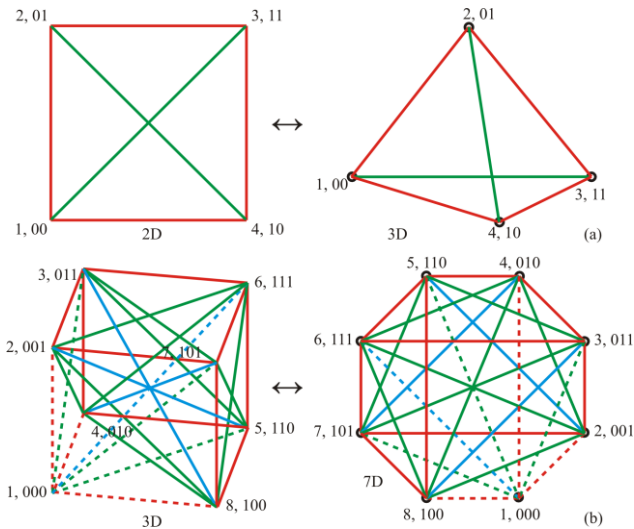

Fig. 2: (a) {2}-cube and 3-simplex; (b) {3}-cube and 7-simplex. Gray encoding. RGB colours denote increasing Hamming distances between addresses of the vertices. Dashed lines symbolise rank 2 compound predicates $A_s$ related with a vertex $m = 1$, as well as the eigenvalues of the cotan Laplacian of [3]-cube.

---

[2] Cf. https://en.wikipedia.org/wiki/Hypercube.
[3] $2^2$ is divisible by 4 and $2^n$ for $n > 2$ is twice divisible by 2 and thus divisible by 4.
[4] http://www.math.uchicago.edu/~may/VIGRE/VIGRE2011/REUPapers/Markov.pdf

Studying properties of this complete graph Pentti Kanerva [PK88] introduced the concept of the Sparse Distributed Memory (SDM), as a mathematical model for memory and learning processes of animals provided with neural networks. Sparseness reflects his hypothesis that not all addresses of this address space are implemented. The main attribute of this model is sensitivity to similarity, meaning that an information can be read back not only by giving the original write address but also by giving one close to it, as measured by the number of mismatched bits (i.e. the Hamming distance between memory addresses).

The model features "knowing that one knows" and "tip of the tongue" phenomena present in biological autonomous learning systems, such as a human brain. According to Kanerva these systems base their operation on an internal model of the world, which they build through experience and a sparse distributed memory is ideal for storing a predictive model of the world [PK88]. Various architectures of artificial neural networks utilising the properties of this graph have been proposed and used for various applications including vision-detecting, robotics, signal detection, etc.

The number of addresses that ere exactly $k$ bits from an arbitrary address $a_m$ is the number of ways to choose $k$ coordinates from a total of $n$ coordinates, and is therefore given by the binomial coefficient ($n : k$). An outstanding property of {$n$}-cube is that the mean Hamming distance between any address $a_m$ and addresses of all vertices (including $m$) is $n/2$ (variation is $n/4$). If an {$n$}-cube was inscribed in the closed $n$-ball most of its vertices would lie at or near the equator. This is called a tendency to orthogonality in SDM language [PK88].

Properties of {$n$}-cube have also been studied by Satosi Watanabe in a framework of "Epistemological Relativity" [SW86]. He noted that each $k^{th}$ bit of an address $a_m$ of a vertex $m$ can be considered a certain Boolean-valued starting predicate $Q_k$. He termed an address $a_m$ a disjoint atomic predicate that can be expressed as a conjunction of starting predicates $Q_k$

$$a_m = \bigcap_{k=1}^{n} a_{mk} Q_k \qquad \text{(ATP)}$$

where $a_{mk} = 1$ if $a_m$ contains $Q_k$ under the conjunction, and $a_{mk} = 0$ if $a_m$ contains the negation of $Q_k$ under the conjunction. In this approach starting predicates $Q_k$ are unrelated to Cartesian coordinates of the vertices. Address $a_7 = 101$ for example can be formed from true predicates $Q_1 = 1$, $Q_2 = 1$, and $Q_3 = 1$ using coefficients $a_{71} = 1$, $a_{72} = 0$, and $a_{73} = 1$.

Watanabe also explored the concept of implicational constraints $Q_k \Rightarrow Q_l$ (which is equivalent to $Q_k = Q_k \cap Q_l$) among the starting predicates. Every time there is a logical constraint of this kind among the $Q_k$'s, one address (ATP) becomes $\emptyset$, and drops out of the list. Any $C$ independent implicational constraints reduces the number of atomic predicates from $M = 2^n$ to [SW86]

$$M = \left(\frac{3}{4}\right)^C 2^n \qquad \text{(ICT)}$$

Watanabe also considered Boolean-valued compound predicates $A_s$ ($s = 1, 2,\ldots, 2^M$), all the logical functions that can be formed from the starting predicates $Q_k$, with connectives of negation, conjunction and disjunction. For $\{n\}$-cube they can also be expressed as a disjunction of vertices

$$A_s = \bigcup_{m=1}^{2^n} A_{sm} \qquad \text{(CPP)}$$

where $A_{sm} = 1$ if $A_s$ contains vertex $m$ under the disjunction and 0 otherwise, so 00…0 denotes the empty set (∅) and 11…1 all $2^n$ vertices (□). Indexing of the predicates $A_s$ is irrelevant in (CPP). Implicational constraints (ICT) make $M = 3^{n/2}$ odd for $C = n/2$ and if $n$ is even. But still $2^M$ remains even.

Compound predicates form a half-ordered Boolean lattice in the sense that there are implicational relations $A_s \to A_t$ between some pairs of the predicates $A_s$, $A_t$. An implicational relation $A_s \to A_t$ is equivalent to $A_{sm} \leq A_{tm}$ for all $m$, in other words, if $A_{sm} = 1$, then $A_{tm} = 1$ for all $m$. Any $A_s$ satisfies $\emptyset \to A_s \to \square$ and the sum of ones in this $2^n$-tuple is the rank $r$ of the predicate $A_s$ or the number of vertices it is built upon[5].

Predicates $A_s$ can thus be thought of as $k$-simplices of a $(2^n-1)$-simplex isomorphic to $\{n\}$-cube with:
(-1)-simplex as the empty set ∅ (rank 0 predicate $A_s$);
0-simplices as vertices (rank 1 or atomic predicates $A_s$);
1-simplices as edges (rank 2 predicates $A_s$);
2- simplices as triangles (rank 3 predicates $A_s$);
3- simplices as tets (rank 4 predicates $A_s$);
and so on up to the single $(2^n-1)$-simplex spanned on all the vertices of the $\{n\}$-cube ($r = k+1$). Nonetheless $k$-simplex contains all its $l$-faces ($0 \leq l \leq k$), while rank $r \geq 1$ predicate $A_s$ contains at least one vertex.

Watanabe [SW86] assumed that any object in the universe satisfies or negates each starting predicate $Q_k$. In other words any object corresponds to a vertex $m$ of $\{n\}$-cube. This assumption is just the identity of indiscernibles ontological principle stating that there cannot be separate objects that have all their properties in common. This principle, however, fails in the quantum domain [SL03]. If this was true a compound predicate $A_s$ of rank $r \geq 2$ (at least an edge) would be shared by $p$ objects if it included $p$ vertices corresponding to these objects. The number of predicates $A_s$ of rank $r \geq 2$ shared by $p$ objects (vertices)

$$N_{r,p} = \binom{2^n - p}{r - p}$$

is the same for any $p$ objects to which the predicates (or simplexes) $A_s$ are applicable (two objects share 1 edge, three objects share 1 face of $(2^n-1)$-simplex, etc.). Watanabe regarded the number of shared predicates as a measure of similarity and the number of not shared predicates as an indication of dissimilarity [SW69]. Therefore any two objects, in so far as they are distinguishable (i.e. correspond to different vertices), are equally similar. This is the Watanabe's famous Ugly Duckling Theorem which, just like quantum theory, defies *common sense*.

As a *corollary* (or rather a relief) to this theorem Watanabe suggested [SW86] that one has to ponderate (give weights to) the predicates $A_s$ to assert the similarity of the objects: for two objects to be similar to each other, they have to share some more important (more weighty) predicates. But it is not convincing: a corollary of a mathematical theorem cannot reduce to some kind of a rule of thumb.

In SDM $\{n\}$-cube models a single neuron, where the vertices play the role of synapses, the points of electric contact between neurons or simply the neuron inputs. The activation function of an artificial neuron defines the output of that neuron given the set of inputs. Only nonlinear activation functions allow neural networks to compute nontrivial problems. Their important characteristic is that they provide a smooth, differentiable transition as input values change, i.e. a small change in input produces a small change in output. A neuron fires when the activation function exceeds a specific activation threshold and is ready for a subsequent perception.

One of the most popular nonlinear activation function is the logistic (aka sigmoid or soft step) one

$$f(x) = \frac{1}{e^{-\mu x} + 1}$$

where $x$ is the weighted sum of the neuron inputs and $\mu$ is a parameter[6]. Rectified linear is another one.

If we now assume that vertices of $\{n\}$-cube are synapses, then it remains to define the activation function and assume a certain activation threshold. We note that

$$N_{r,1} = \binom{2^n - 1}{r - 1} \text{ out of } \binom{2^n}{r}$$

rank $r$ available compound predicates $A_s$ (CPP) of rank $r \geq 1$ are related with one vertex $m$. For $\{3\}$-cube, for example, one vertex is related with 1/8 rank 1, 7/28 rank 2, 21/56 rank 3, 35/70 rank 4, 35/56 rank 5, 21/28 rank 6, 7/8 rank 7 predicates, and with the whole set of vertices of $\{3\}$-cube (cf. Fig. 2(b)). This linear sequence can be used to define rectified linear activation function of $\{3\}$-cube neuron associated with this vertex $m$. In general

$$N_{r,p} = \sum_{l=1}^{p} \binom{2^n - l}{r - 1} \text{ out of } \binom{2^n}{r}$$

rank $1 \leq r \leq 2^n$ available predicates $A_s$ is related with $p$ vertices, where binomial coefficient is defined in terms of a falling factorial for $r-1 > 2^n-l$ (i.e. zero), as illustrated in Fig. 3 showing sequences

$$f(r,p) = \sum_{l=1}^{p} \binom{2^n - l}{r - 1} \bigg/ \binom{2^n}{r} \qquad \text{(CAF)}$$

---

[5] To note in passing, a compound predicate $A_s$ of rank r ≥ 2 related with a vertex $m$ is in general not the same as a circle $O(r, a_m)$ with radius $r$ and center $a_m$ according to SDM definition. Only a circle with radius $n$ equals rank $n$ predicate $A_s$ and contains the set of all $2^n$ vertices. A circle with radius zero is $a_m$ not the empty set.

[6] The values of interest for $\mu$ seem to be those in the interval [0,4] (or [-2,4]). In parts of this range the logistic map $x_{n+1} = \mu x_n (1 - x_n)$ which is analogous to the logistic function derivative $f'(x) = \mu f(x)(1 - f(x))$ displays intermittent (irregular alternation of periodic and chaotic dynamics) behaviour.

for {3}-cube and {4}-cube.

We see that (CAF) is undefined if no vertex is active ($r = 0$), linear for $p = 1$, and for $p > 1$ is nonlinear and similar to the logistic activation function up to $p = 2^n$, where it becomes unit constant function. These sequences are the same for any $p$ vertices but remain specific for vertices that are comprised in the set of size $|p|$. In this way (CAF) can be thought of as the neuron's vertex dependent activation function parametrised by $p \geq 1$.

Extensions of this model to artificial neural networks of {$n$}-cubes should certainly take the implicational constraints (ICT) into account. It may well be that they are formed in a biological neural network of {$n$}-cubes in a result of semiosis reducing the number of available vertices. This conjecture is supported by the fact that interpreting information requires some predefined language (or structure) that could certainly be provided by implicational constraints, be it the language of hormones, pheromones, barking, or English.

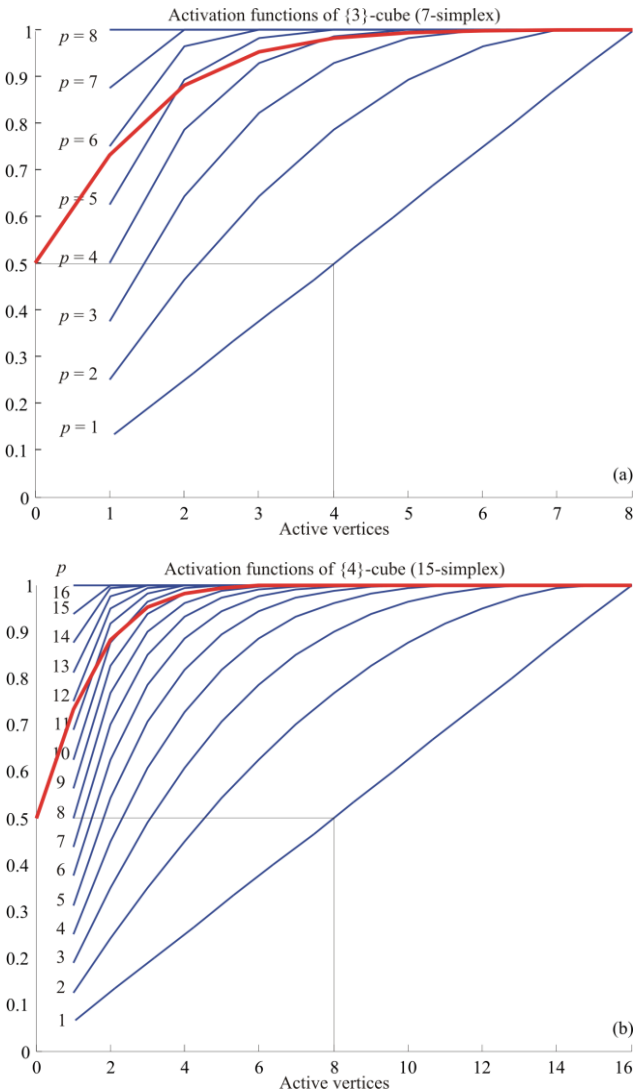

Fig. 3: Activation functions of {$n$}-cube and the logistic activation function ($\mu = 1$). (a) {3}-cube, (b) {4}-cube.

A neuron as a living biological cell is a dissipative structure. That hints another cube.

## 3.  [$n$]-cube

2-dimensional triangulated surface allows to define a discrete cotan-Laplace operator assuming that every relation (edge) between vertices $i$ and $j$ carries a real-valued weight

$$\omega_{ij} = \frac{1}{2}\left(\cot \alpha_{ij} + \cot \beta_{ij}\right) \qquad \text{(OLP)}$$

where $\alpha_{ij}$ and $\beta_{ij}$ are the angles opposite the edge between vertices $i$ and $j$. This equation, called cotangent formula, has been derived in many different ways, and rediscovered many times[7] over the years [KC19a].

If the vertices are ordered, this ordering can be used to induce the orientation of the edges and angles $\alpha_{ij}$ and $\beta_{ij}$. as $i < j \Leftrightarrow i \rightarrow j$. Then if $\alpha_{ij}$ and $\beta_{ij}$ are directed towards vertex $j$ weights, $\omega_{ij}$ are positive (OLP), while if they are directed towards vertex $i$ weights, $\omega_{ij}$ are negative

$$\omega_{ij} = \frac{1}{2}\left(\cot(-\alpha_{ij}) + \cot(-\beta_{ij})\right) = \\ = -\frac{1}{2}\left(\cot(\alpha_{ij}) + \cot(\beta_{ij})\right) \qquad \text{(OLN)}$$

For a "locally disk-like" triangulated manifold allowing at most two triangles incident to an edge, a discrete cotan-Laplace operator acting on a function $u: V \rightarrow \mathbb{R}$, where $V$ is the a vertex set of this graph is defined [MW17] as

$$(Lu)_i = \sum_{j \sim i} \omega_{ij}\left(u_i - u_j\right) \qquad \text{(WDL)}$$

where the sum ranges over all vertices $j$ that are related with the vertex $i$. This allows for representing the linear operator $L$ as a matrix

$$L_{ij} = \begin{cases} -\omega_{ij} & \text{if } i \neq j \text{ and they are related} \\ \sum_{l \sim i} \omega_{il} & \text{if } i = j \\ 0 & \text{if } i \text{ and } j \text{ are unrelated} \end{cases} \qquad \text{(WDM)}$$

which is called the weakly defined discrete Laplace matrix or just cotan Laplacian [MW17][8]. Zeroing $\omega_{ij}$ for unrelated vertices encapsulates the locality of action of the Laplacian operator (WDL): changing the value of $u_j$ at a vertex $j$ does not alter the value $(Lu)_i$ at vertex $i$ if these vertices are unrelated.

Noting in passing: in a distance matrix $D$ (weighted adjacency matrix), such as (DMT), entries pertaining to unrelated vertices should be set to infinity or a suitable large value, as zero in these locations would be incorrectly interpreted as an edge with no distance, cost, etc. Both (OLP) and (OLN) guarantee setting such entries to infinity if $\alpha_{ij}$ and $\beta_{ij}$ are 0 or $\pi$. One could say that in this case such relation is *zero* relation or *flat* relation.

---

[7] The same physical laws can be derived from numerous starting points and numerous assumptions. Take the entropic gravity derivation of the law of gravitation as an example.

[8] This is in general not the same as Kirchhoff matrix ($G - E$), where $G$ is the (diagonal) degree matrix and $E$ the adjacency matrix of the graph.

One may define the gradient $\nabla u_{ij}$ between the related vertices $i$ and $j$ as the finite difference $(u_i - u_j)$. Accordingly, one defines discrete Dirichlet energy of $u$ as

$$E_D[u] \doteq \frac{1}{2}\sum_{e \in E} \omega_{ij}(u_i - u_j)^2 = \frac{1}{2}u^T L u \quad \text{(DEN)}$$

where the sum ranges over all relations between the vertices. Solving discrete cotan-Laplace equation $Lu = 0$ for all vertices and subject to appropriate boundary conditions, is equivalent to solving the variational problem of finding a function $u$ that satisfies the boundary conditions and has minimal Dirichlet energy (DEN).

The Delaunay triangle mesh (dual to the Voronoi one) features a number of optimality properties: the triangles are the "fattest" possible [RCEA09], it maximizes the minimal angles in the triangulation, and more importantly the Delaunay triangulation of the set of vertices of $n$-manifold minimizes the Dirichlet energy of any piecewise linear function $u$ over this point set (Rippa's theorem).

Cotan Laplacian $L$ is (the list is not exhaustive):
- singular (non-invertible), as it has zero eigenvalue;
- symmetric (self-adjoint) ($\omega_{ij} = \omega_{ji}$) and thus has real eigenvalues and orthogonal eigenvectors;
- $L$ with positive weights (OLP) is always positive semi-definite, and $L$ with negative weights (OLN) is always negative semi-definite (these are not necessary conditions however [MW17]);
- has only constant functions $u$ in its kernel [MW17];
- satisfies the maximum principle since

$$\frac{1}{L_{ii}}\sum_{l \sim i}\omega_{il} = \sum_{l \sim i}\omega_{il} \bigg/ \sum_{l \sim i}\omega_{il} = 1$$

and thus $u_i$ is a convex combination of its related neighbours $u_j$ for discrete harmonic functions;
- for any $2^n \times (2^n-1)$ matrix $B$, one has

$$\det(B^\dagger L B) = p_L'(0)\left|\det(B\ x_i)\right|^2$$

where $(B\ x_i)$ denotes the $2^n \times 2^n$ matrix with right column containing an $L$ eigenvector $x_i$ and all remaining columns given by $B$ and $p_L'(0)$ is the derivative of the characteristic polynomial of $L$ at its zero eigenvalue [PDEA20][9], the product of all the remaining non-zero eigenvalues of $L$;
- the spectrum of the cotan Laplacian obtains its minimum on a Delaunay triangulation in the sense that the $i$-th eigenvalue of the cotan Laplacian of any triangulation of a fixed point set in the plane is bounded below by the $i$-th eigenvalue resulting from the cotan Laplacian associated with the Delaunay triangulation of the given point set [RCEA09].

The empty circle property of Delaunay triangulation implies that an interior edge is Delaunay edge iff $\alpha_{ij} + \beta_{ij} \leq \pi$, which is equivalent to $\sin(\alpha_{ij} + \beta_{ij}) \geq 0$ or $\cot(\alpha_{ij}) + \cot(\beta_{ij}) \geq 0$. It has been shown in [CL72] that any convex quadrilateral formed by two adjacent triangles which does not satisfy the empty circle property may be made Delaunay by "flipping" the diagonal edge of the quadrilateral, common to the two triangles, to the opposite diagonal. This is called a Delaunay flip and a sequence of Delaunay flipping will always converge to a Delaunay triangulation [RCEA09]. The borderline case for a Delaunay flip is obviously a rectangle, or a square in particular, having both Delaunay diagonal edges. This leads to the following observation (cf. Fig. 4).

---
[9] Lemma 13 (Cauchy-Binet type formula).

**Theorem 1**: Any one of the two possible triangulations on 2-faces of $n$-cube produces the same cotan Laplacian.

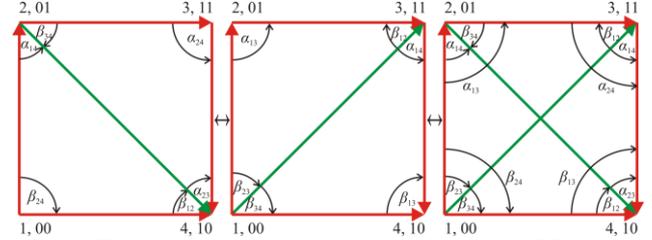

Fig. 4: Two diamond graphs on 2-faces of $n$-cube (OLP) have the same cotan Laplacian that defines [$n$]-cube cotan Laplacian.

**Proof 1**: It is easy to see. For $n > 2$ there are at least two angles opposite to each 2-face edge. For 2-face edges $\omega_{ij} = 1$ since $\cot(\pi/4) = 1$ for each angle, regardless of the arrangement of 2-face diagonals, and for 2-faces diagonals $\omega_{ij} = 0$ since both angles opposite to a diagonal are right angles ($\cot(\pi/2) = 0$). Therefore the coefficients $L_{ij}$ of such matrix are zero not only for unrelated pairs of vertices $i$ and $j$ but for all pairs of vertices, except for 2-face edges, where $L_{ij} = -1$ and diagonal coefficients where $L_{ii} = n$, as any vertex of $n$-cube is connected with $n$ 2-face edges (for (OLN) signs are reversed). □

Technically cotan Laplacian for [$n$]-cube can be simply produced from $n$-cube Hamming distances matrix $D$ (DMT) by setting $L_{ii} = n$, negating ones and zeroing all the other entries (for (OLN) signs are reversed).

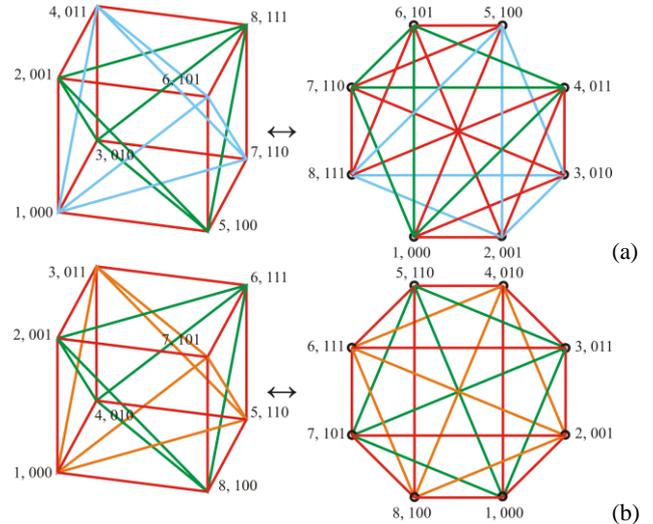

Fig. 5 [3]-cube. (a) binary, encoded (b) Gray encoded. For Gray encoding colors represent even (orange) and odd (green) arrangements of 2-face diagonals.

**Definition 1**: [$n$]-cube shown in Fig. 5 is defined as $n$-cube with triangulated 2-faces. From Theorem 1 it follows that even or odd arrangement of 2-face diagonals is irrelevant for the discrete cotan Laplacian, so any, both or none of these triangulations can be equivalently used. It is equivalent to non-triangulated $n$-cube. Since $\cot(\pi/2) = 0$ only the ordering of the vertices (Gray, binary, random, etc.) changes the form of the cotan Laplacian. Binary encoding can obviously be left-msb, right-msb, or any other sequence out of $n!$ possible orderings of bits.

There is a formal problem with this definition of [$n$]-cube, as ($n$-1) 2-faces adjoin any given boundary edge. So for $n \neq 3$ there are more than two triangles per boundary edge, and it is undefined which $\omega_{ij}$ should one take to (WDM) even though $\alpha_{ij}$ always equals $\pi/4$. To overcome this we could demand the graph of [$n$]-cube to be Eulerian circuitable, which guarantees that each edge would be traversed only once. Degree of a vertex of [$n$]-cube (if both triangulations are present) is $(2|e|+4|f|))|v|$ or ($n$+1 : 2) (triangular numbers) and by Euler's Theorem this graph has an Eulerian circuit if and only if this is even, which holds only if $\lfloor(n+1)/2\rfloor$ is even, that is, e.g., for $n = 3$ or $n = 4$ (cf. Fig. 3). For Gray encoded [3]-cube an Eulerian circuit can be formed on 24 edges between vertices 362712876815846543142357 3. This swap of parity with every two dimensions is a property also present in (VNB) and (SNB) formulas.

Also [1]-cube and [2]-cube are exceptions; [2]-cube is the only triangulation with a boundary, while [1]-cube has no triangulation at all. But even in these degenerate cases the cotan Laplacian can be obtained by modification of the corresponding matrix $D$ (DMT).

**Theorem 2**: For $n > 2$ The eigenvalues of the cotan-Laplacian of any one of two triangle meshes on 2-faces of [$n$]-cube correspond to twice the binomial distribution of Hamming distances between any address $a_m$ and addresses of all vertices (including $m$) as in {$n$}-cube. E.g. {0, 2, 2, 2, 4, 4, 4, 6} for {3}-cube, etc. For (OLP) 0 is the minimum eigenvalue; for (OLN) 0 is the maximum eigenvalue and the distances are negative. For $n = 2$ the eigenvalues are Hamming distances {0, 1, 1, 2} unless we assume that $L$ is just a form of $D$ for $n \leq 2$ (or equivalently that [2]-cube has no boundary). For $n \geq 2$ the spectrum is degenerate (eigenvalues of $L$ are not distinct).

**Proof 2**: Direct calculation shows that it is true. □

Theorem 2 interesting since (WDL) does not require neither the notion of a metric nor the vertices coordinates. The binomial distribution of the Hamming distances between the coordinates of any vertex $m$ of {$n$}-cube and all the other vertices arises naturally from the angles of 2-face triangulation on [$n$]-cube.

**Theorem 3**: The maximum of the absolute values of all nontrivial eigenvalues of adjacency matrix of regular [$n$]-cube is

$$\lambda_{\max}(n) = \max_{i \neq 1}|\lambda_i| = \sum_{k=1}^{n} k - 2n = -\frac{n(3-n)}{2}$$

(negation of A080956 OEIS sequence). Thus ($n + (n : 2)$)-regular [$n$]-cube is Ramanujan graph for n < 6, that is[10]

$$-\frac{n(3-n)}{2} \leq 2\sqrt{n + \binom{n}{2} - 1} \quad \text{iff} \quad n < 6 \quad \text{(RMG)}$$

**Proof 3**: Direct calculation shows that it is true. □

**Theorem 4**: Integer values of the Ramanujan graph constant for [$n$]-cube (RHS of (RMG)) form A075848 OEIS sequence given by

$$\frac{3}{2\sqrt{2}}\left[\left(3+2\sqrt{2}\right)^k - \left(3-2\sqrt{2}\right)^k\right]$$

and $n$'s yielding these integer values form A072221 OEIS sequence given by

$$\frac{3}{4}\left[\left(3+2\sqrt{2}\right)^k + \left(3-2\sqrt{2}\right)^k\right] - \frac{1}{2}$$

where $k \in \mathbb{N}_0$.

**Proof 4**: Direct calculation shows that it is true. □

The cotan Laplacian $L$ of [$n$]-cube has the following additional properties to the cotan Laplacian (WDM) of any 2-dimensional triangulated surface (the list is not exhaustive; certain properties of $L$ are possibly duplicated; (OLP) is assumed and binary or Gray encoding is necessary, unless indicated otherwise):

- is bisymmetric, that is both symmetric ($L = L^T$) and centrosymmetric ($LJ = JL$, where $J$ is the exchange matrix having $2^{n-1}$ eigenvalues +1 and $2^{n-1}$ eigenvalues -1);
- has $2^{n-1}$ symmetric orthonormal eigenvectors $x$ such that $Jx = x$ [CB76];
- has $2^{n-1}$ antisymmetric (skew-symmetric) orthonormal eigenvectors $x$ such that $Jx = -x$ [CB76];
- is an $M$-matrix (real $Z$-matrix);
- its eigengap (the difference between two successive eigenvalues) equals 2;
- its spectral radius equals $2n$;
- its trace equals $n2^n$ ($n2^{n-1}$ is the number of edges of $n$-cube);
- $n$-cube is $n$-regular, so the normalized Laplacian $\mathcal{L} = L/n = I - A/n$ (but it doesn't have such an interesting spectrum of integers);
- equals Kirchhoff Laplacian

$$L = nI - E$$

where $E$ the adjacency matrix of $n$-cube that can be written as

$$E = \begin{bmatrix} B & C \\ C & B \end{bmatrix}$$

where $B$ is a bisymmetric hollow matrix, $C = I$ for binary encoding and $C = J$ for Gray encoding, and $E$ is self-similar, that is $L$ is built from $2^{n-1}$ 2×2 matrices

$$L_{2\times 2} = \begin{bmatrix} n & -1 \\ -1 & n \end{bmatrix}$$

---

[10] The maxima of the absolute values of all nontrivial eigenvalues of adjacency matrices of regular $n$-cubes form A023444 ($n$-2) OEIS sequence, so $n < 7$; for {$n$}-cubes they are -1, so all {$n$}-cubes are Ramanujan graphs.

- (having eigenvalues $\lambda(L_{2\times 2}) = n \pm 1$) on diagonal and corresponding number[11] of $C$ matrices;
- is orthogonally similar to the matrix $O$ [CB76]

$$O = \begin{bmatrix} A - JC & 0 \\ 0 & A + JC \end{bmatrix}$$

that is $L = K^T O K$ with

$$K = \frac{1}{\sqrt{2}} \begin{bmatrix} I & -J \\ I & J \end{bmatrix}$$

and $\det(KLK^T) = \det(A - JC)\det(A + JC)$, where $A - JC$ is linearly independent and $A + JC$ is linearly dependent and eigenvalues of $A - JC$ correspond to antisymmetric eigenvectors of $L$ and eigenvalues of $A + JC$ correspond to symmetric eigenvectors of $L$ [CB76], as similarity transformations preserve eigenvalues;
- eigenvectors are defined up to a phase, that is, if $Lx = \lambda x$ then $e^{i\theta}x$ is also an eigenvector of $L$, and specifically so is $-x$ (where $\theta = \pi$); this somehow introduces quantum theory to [n]-cube picture;
- for $k \ne 0$ and $k \ne n$ the eigenvalues $\lambda_{(k)}$ are degenerate, so the eigenvectors corresponding to these eigenvalues have an additional freedom of rotation, that is to say any linear (orthonormal) combination of eigenvectors sharing such eigenvalue $\lambda_{(k)}$ (in the degenerate subspace), are themselves eigenvectors (in the subspace);
- columns and rows of $L$ are linearly dependent, so that Laplace equation for [n]-cube, which is a homogeneous system

$$Lu = 0 \qquad \text{(LEQ)}$$

has trivial solution;
- has only constant vectors in its kernel, that is the other solutions of (LEQ) are constant vectors $u_v$ in any encoding; therefore if $u_p$ is any specific solution to the Poisson's equation for [n]-cube, which is just the linear system

$$Lu = f \qquad \text{(PEQ)}$$

then the entire solution set can be described as $\{u_p + u_v\}$, where $u_v$ is a constant vector solving (LEQ);
- it is a linear map (additivity and scalar multiplication);
- Gray encoded $L$ can be Cholesky decomposed;
- binary encoded $L$ can be LDL decomposed;

Solving (PEQ) for binary encoded [3]-cube for $f_i = |1|$ with $2^{n-1}$ of $f_i = 1$ and $2^{n-1}$ of $f_i = -1$ ((8 : 4) = 70 possibilities) yields 2/3 as the minimum norm of potential vector $u$ for ones at vertices 2, 3, 5, 8 and minus ones at vertices 1, 4, 6, 7 or vice versa, as shown in Fig. 5 (a). This hints yet another cube.

### III. Fourth cube

**Definition 2**: $2^n$-cube is defined by $2^n$ $n$-cubes sharing common origin $(00...0)$. This structure has $3^n$ vertices and

$$\binom{n}{k} 3^{n-k} 2^k$$

---
[11] $2^{n-k} C(2^k \times 2^k)$, $k = 1,...,n-1$.

$k$-faces, $k = 0, 1, ..., n$ (A038220 OEIS sequence) and sum of all $k$-faces is $5^n$. The number of addresses having the same $k$-norm (i.e. the same $k$ ones) is

$$\binom{n}{k} 2^k$$

which is A013609 OEIS sequence.

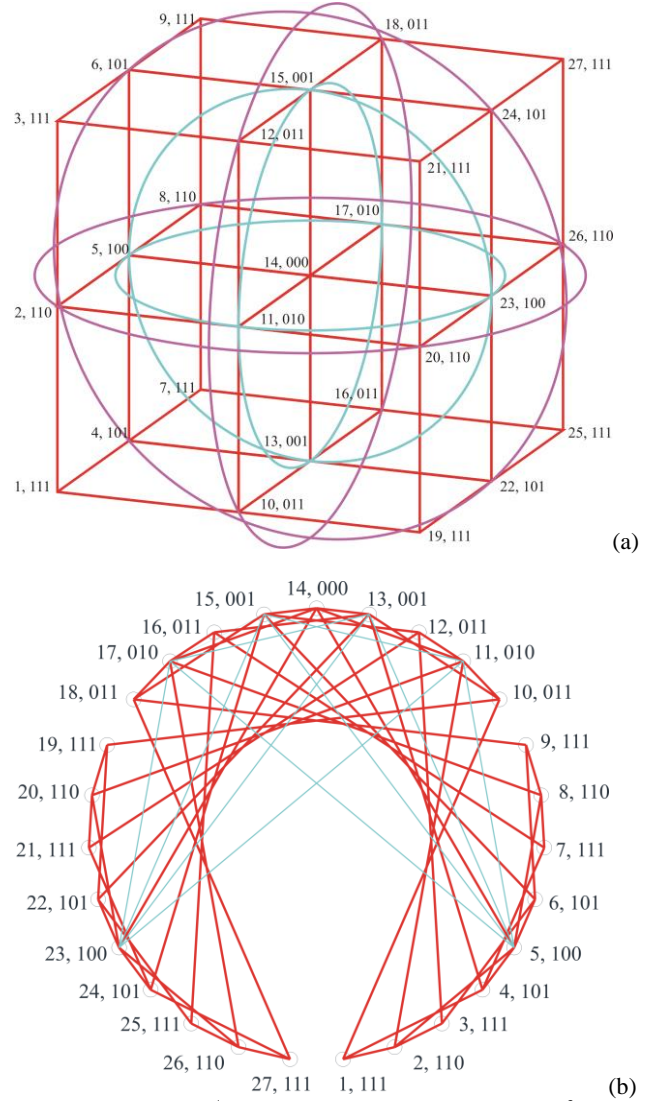

Fig. 6: (a) unit and $\sqrt{2}$ radii (3)-balls in binary encoded $2^3$-cube, (b) graph of this structure.

**Definition 3**: $[2^n]$-cube is defined analogously as $2^n$-cube, while maintaining the properties of [n]-cube: it contains only boundary edges and 2-face diagonals carrying weights given by (WDM).

**Definition 4**: $\{2^n\}$-cube is defined analogously as $2^n$-cube, while maintaining the properties of $\{n\}$-cube: edges carry Hamming distances between addresses of vertices.

These structures resemble Cartesian coordinate system of $\mathbb{R}^n$, as shown in Fig. 6, wherein the signs of coordinates are provided by indexation of the vertices. In other words $2^n$-cube provides a bijective relation between indices of $3^n$ vertices it comprises and their coordinates in $n$-dimensional space. $2^n$-cube corresponds to unit $n$-cube, as shown in Fig. 9 for $n = 3$; there is only one central vertex 14 having address $a_{14} = [0,0,0]$ from which one can recur-sively walk to six 1-norm vertices 13, 15; 11, 17 and 5, 23 by adding $3^k$ to or subtracting $3^k$ from the vertex index with $k = 0,1,2$ to and modifying the corresponding bit of the new vertex address, and so on to 2-norm vertices up to 3-norm vertices 1,3;7,9;19,21;25,27 having address [1,1,1].

Nonetheless, the Hamming distance between vertices having different indices but the same addresses is zero, even though Euclidean metric between the points having Cartesian coordinates induced by these vertices is between 2 to $2\sqrt{n}$. The Hamming distance between addresses of vertices 1 and 27 of $2^3$-cube, for example, is 0, even though Euclidean metric between Cartesian coordinates induced by these vertices is $2\sqrt{3}$. This provides interesting properties of the distance matrix of $\{2^n\}$-cube. Distance matrix of $\{2^2\}$-cube, for example, defines a torus.

Certain properties of these graphs are discusses below.

1. **$2^n$-cube**

**Theorem 5**: Spectrum of eigenvalues of adjacency matrix of $2^n$-cube in Gray or binary encoding is symmetric, degenerate (for $n \geq 2$), and includes all multiplicities of $\sqrt{2}$ from $-n\sqrt{2}$ to $n\sqrt{2}$. Multiplicities of the same eigenvalues form trinomial triangle.

**Proof 5**: $2^n$-cube is bipartite, so the spectrum is symmetric. Direct calculation shows that it has the remaining properties. □

2. **$[2^n]$-cube**

**Theorem 6**: Cotan Laplacian (WDM) spectrum of eigenvalues of $[2^n]$-cube in Gray or binary encoding is degenerate and includes all integers from 0 to $3n$ without the eigenvalue of $3n$-1. Multiplicities of the same eigenvalues form A038717 OEIS sequence.

**Proof 6**: Direct calculation shows that it is true. □

Some of the further properties of the cotan Laplacian of this graph are (again (OLP) is assumed and binary or Gray ordering of vertices is necessary to center the origin):
- it is bisymmetric;
- it is odd dimensional ($3^n \times 3^n$);
- its diagonal entries span from $n$ to $2n$;
- its spectral gap (the difference between two largest eigenvalues) equals 2; otherwise the eigengap is 1;
- its spectral radius equals $3n$;
- it has $\lceil 3^n/2 \rceil$ symmetric orthonormal eigenvectors [CB76];
- it has $\lfloor 3^n/2 \rfloor$ antisymmetric orthonormal eigenvectors [CB76];
- it can be written as

$$L = \begin{bmatrix} A & x & C^T \\ x^T & 2n & x^T J \\ C & Jx & JAJ \end{bmatrix}$$

where $A = A^T$, and $x$ is a column vector having $\lceil 3^n/2 \rceil$ rows with -1 in $\lceil 3^n/2 \rceil + 1 - 3^k$, where $k = 0,1,\ldots,n-1$, row's entries and zeros elsewhere; -1's in the central row and column of $L$ correspond to kissing vertices of $n$-ball inscribed inside $[2^n]$-cube;

- it is orthogonally similar to the matrix $O$ [CB76]

$$O = \begin{bmatrix} A - JC & 0 & 0 \\ 0 & 2n & \sqrt{2}x^T \\ 0 & \sqrt{2}x & A + JC \end{bmatrix}$$

that is $L = K^T O K$ with

$$K = \frac{1}{\sqrt{2}} \begin{bmatrix} I & 0 & -J \\ 0 & \sqrt{2} & 0 \\ I & 0 & J \end{bmatrix}$$

- it equals Kirchhoff Laplacian $L = G - E$, where $G$ is the (diagonal) degree matrix and $E$ the adjacency matrix of $2^n$-cube that can be written as

$$E = \begin{bmatrix} B & C & 0 \\ C & B & C \\ 0 & C & B \end{bmatrix}$$

where $B$ is a bisymmetric hollow matrix, $C = I$ for binary encoding and $C = J$ for Gray encoding, and $E$ is self-similar, that is $L$ is built from $3^{n-1}$ $3\times 3$ matrices

$$L_{3\times 3} = \begin{bmatrix} k-1 & -1 & 0 \\ -1 & k & -1 \\ 0 & -1 & k-1 \end{bmatrix}$$

(having eigenvalues $\lambda(L_{3\times 3}) \in \{k-2, k-1, k+1\}$) on diagonal, where $k = n+1, n+2,\ldots,2n$ for $n \geq 2$ and corresponding number[12] of $C$ matrices.

$[2^n]$-cube enables to inscribe unit $n$-ball inside it[13], having surface kissing all vertices distanced one bit from the centre, as shown in Fig. 6(a). In fact it enables to inscribe all $n$-balls having $\sqrt{k}$ radii, where $k = 1,2,\ldots,n$, having surfaces kissing all vertices distanced $k$ bits from the centre. $k = \{1, 2\}$ are of particular importance as they pertain to triangular surface meshes yielding simplified form of the cotan Laplacian for $\alpha_{ij} = \beta_{ij} = \pi/2$. This can be seen in $[2^2]$-cube as shown in Fig. 7 that is closely related both to the *quadratura circuli* and the lune of Hippocrates ancient problems.

In [VNB08] the construction shown in Fig. 7 is used to interpret the proposed mathematical definition of the fine structure constant

$$\alpha^{-1}(\pi) = 4\pi^3 + \pi^2 + \pi^1 = 137.036\ 303\ 776\ldots$$

close to CODATA (2018) recommended value of

$$\alpha^{-1} = \frac{4\pi\varepsilon_0 \hbar c}{e^2} = 137.035\ 999\ 084\ldots$$

According to this interpretation the tetragonal substitution of $\pi$ can be only the $\pi \sim 4$ as outside the square-measure

---

[12] $4 \cdot 3^{n-k} C(3^{k-1} \times 3^{k-1})$, $k = 2,3,\ldots,n$.
[13] Disc in 4 squares, ball in 8 cubes, etc. Unit $n$-ball volume and unit $n$-ball surface attain their maxima respectively in $n = 5$ and $n = 7$ dimensions ($2n$-1 = 5 and $2^n$-1 = 7 for $n = 3$).

and π~2 as the inside square-measure of the "generative circle" with unit radius

$$\frac{\alpha^{-1}(4)}{2} - 1 = 137 < \alpha^{-1}(\pi) < 138 = \frac{\alpha^{-1}(4)}{2}$$

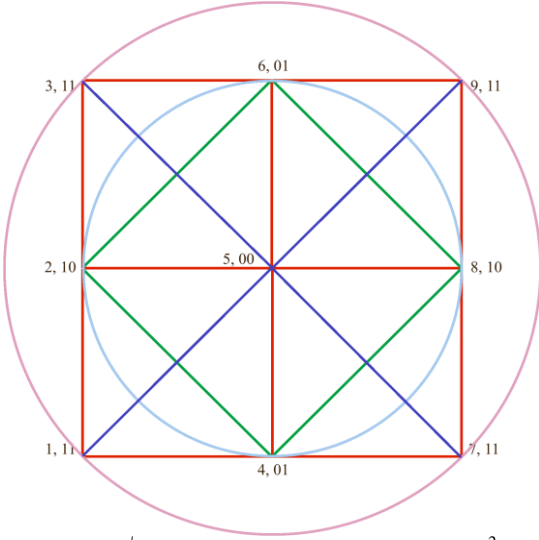

Fig. 7: Unit and √2 radii (2)-balls in binary encoded $[2^2]$-cube or $\{2^2\}$-cube. Distance matrix of $\{2^2\}$-cube defines a torus.

### 3. $\{2^n\}$-cube

**Theorem 7**: Distance matrix for binary and Gray encoded $\{2^n\}$-cube is the same.

**Proof 7**: Direct calculation shows that it is true. Both binary and Gray encoding give the same norms of the addresses with $\lfloor 3n/2 \rfloor + 1$ as the origin. □

Spectrum of eigenvalues of the bisymmetric distance matrix of $\{2^n\}$-cube is degenerate for $n \geq 2$ and has distinct irrational minimum negative eigenvalue and distinct irrational maximum positive eigenvalue given by

$$\lambda_{\min}(n) = \left[2(n-1) - \sqrt{2(2n+1)(n+2)}\right]3^{n-2}$$

$$\lambda_{\max}(n) = \left[2(n-1) + \sqrt{2(2n+1)(n+2)}\right]3^{n-2}$$

as illustrated in Fig. 8.

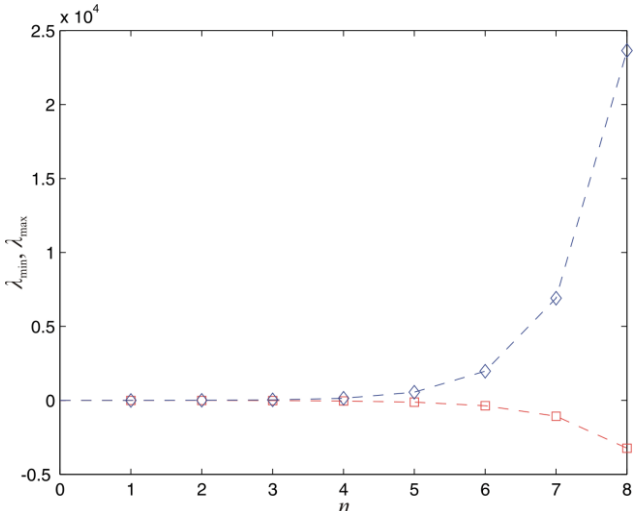

Fig. 8: $\lambda_{min}$ and $\lambda_{max}$ eigenvalues of distance matrix of $\{2^n\}$-cube.

The sums of these eigenvalues

$$\lambda_{\min}(n) + \lambda_{\max}(n) = 4(n-1)3^{n-2}$$

form OEIS integer sequence A120908, while the products

$$\lambda_{\min}(n)\lambda_{\max}(n) = -2n3^{2n-2}$$

form integer sequence which is close to the OEIS sequence A288834 (every 2nd entry's modulus agree). Also

$$\left|\lambda_{\min}(7) + \lambda_{\max}(7)\right| = \left|-\lambda_{\min}(4)\lambda_{\max}(4)\right| = 5832$$

is the only common value between those sums and products.

Furthermore the spectrum contains $n$-1 integer negative eigenvalues given by

$$\lambda_2 = -4 \cdot 3^{n-2}$$

(which is the opposite of the OEIS sequence A003946 for $n \geq 2$) and $3^n$-$n$-1 (this is OEIS sequence A060188) zero eigenvalues. Sum of all distance matrix eigenvalues is zero (obviously product also).

## IV. Fourth cube vs. Vector Equilibrium

"The teacher stood at the blackboard, made a little dot, and said, «This is a point; it doesn't exist». (So she wiped that out.) Then she drew a whole string of them and called it a 'line'. Having no thickness, it couldn't exist either. Next she made a raft out of these lines and came up with a 'plane'. I'm sorry to say it didn't exist either, sighs Bucky. She then stacked them together and got a 'cube', and suddenly that existed. I couldn't believe it; how did she get existence out of nonexistence to the fourth power? So I asked. «How old is it? » She said. «Don't be naughty.»…It was an absolute ghost cube."

Buckminster Fuller [AE87] (p. 10)

The dimension $n$ of a space is a natural number, the minimum number of coordinates defining a point of that space. $n$ = -1 is a dimension of the void, the empty set ∅ or (-1)-simplex. If one assumes that the void contains a single point ($n$ = 0) it is no longer the void. The 1st Peano axiom for the natural numbers (0 is the 1st natural number) results in that $n$ > -1 and induces all the remaining dimensions closed under a successor function $S(n) := n+1$. Therefore primordial Big Bang singularity (a point) does not imply an expansion of space and time but an expansion of dimensions.

An $n$-dimensional simplicial companion of $2^n$-cube is $n$-vector equilibrium, a structure featuring radial equilateral symmetry (circumradius equals the edge length) for any $n$ defined by Buckminster Fuller for $n$ = 3. External vertices of $n$-vector equilibriums are vertices of Stott expanded regular simplices by its dual ones, polytopes forming a dimensional sequence listed in Table 1[14]. All these polytopes (for $n \geq 1$) are 3 layered stacks of vertex layers. The convex hull of the 1st layer is the corresponding ($n$-1)-simplex, which thus is a true facet of that polytope, the convex hull (cross-section) of the 2nd, equatorial, layer is that expanded ($n$-1)-simplex, and the convex hull

---
[14] http://hi.gher.space/forum/viewtopic.php?f=32&t=2207

of the 3$^{rd}$ layer is the ($n$-1)-simplex, dual to the ($n$-1)-simplex of the 1$^{st}$ layer.

| $n$ | Dynkin diagram | name | $|v|$ |
|---|---|---|---|
| -1 | ∅ | empty set | 0 |
| 0 | 1 | point | 0 |
| 1 | 3 | two line segments | 2 |
| 2 | x3x={6} | regular hexagon | 6 |
| 3 | x3o3x=co | cuboctahedron | 12 |
| 4 | x3o3o3x=spid | runcinated 4-simplex | 20 |
| 5 | x3o3o3o3x=scad | stericated 5-simplex | 30 |
| 6 | x3o3o3o3o3x=staf | pentellated 6-simplex | 42 |
| 7 | x3o3o3o3o3o3x=suph | hexicated 7-simplex | 56 |
| 8 | x3o3o3o3o3o3o3x=soxeb | heptellated 8-simplex | 72 |

Table 1: Stott expanded regular simplices by its dual ones.

The number of the external vertices (other than the origin) of the $n$-vector equilibrium is given by

$$|v(n)| = n(n+1)$$

This is A279019 OEIS sequence; least possible number of diagonals of simple convex polyhedron with n faces. This is also twice the sum of the natural numbers less than or equal $n$

$$|v(n)| = 2\sum_{k=0}^{n} k$$

It can also be obtained by recurrence relation

$$|v(n)| = |v(n-1)| + 2n$$

with $|v(-1)| := 0$.

For $1 \leq n \leq 3$ $|v(n)|$ equals the kissing number of $\mathbb{R}^n$; the greatest number of non-overlapping unit spheres that can be arranged in $\mathbb{R}^n$ such that they each touch a common unit sphere. But $|v(4)| = 20$, while the kissing number of $\mathbb{R}^4$ is 24. Known lower bounds (and known exact values for $n = 8$ and $n = 24$) on the kissing number for $n > 4$ are also larger than $|v(n)|$.

Comparison of $2^n$-cubes and $n$-vector equilibrium is illustrated in Table 2 for $n \leq 3$.

| $n$ | $2^n$-cube | $n$-vector equilibrium |
|---|---|---|
| -1 | | ∅ |
| 0 | |  |
| 1 | 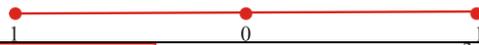 | 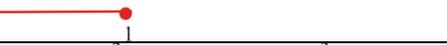 |
| 2 | 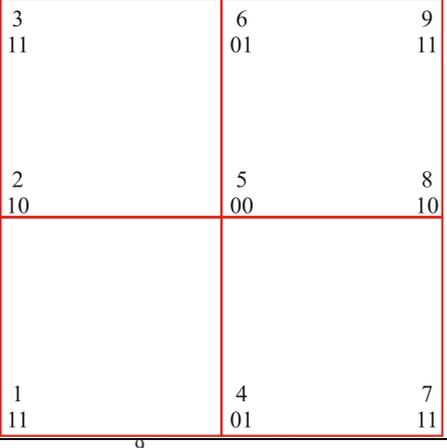 | 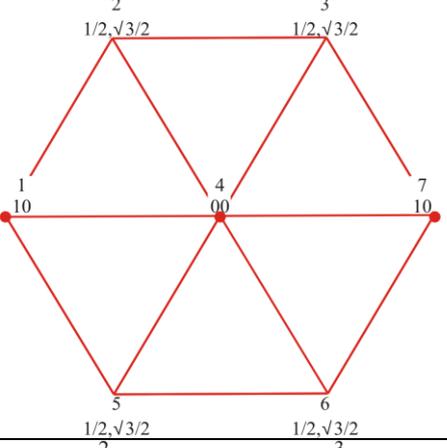 |
| 3 | 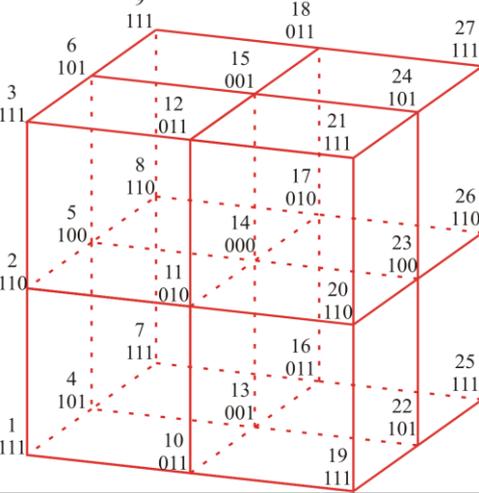 | 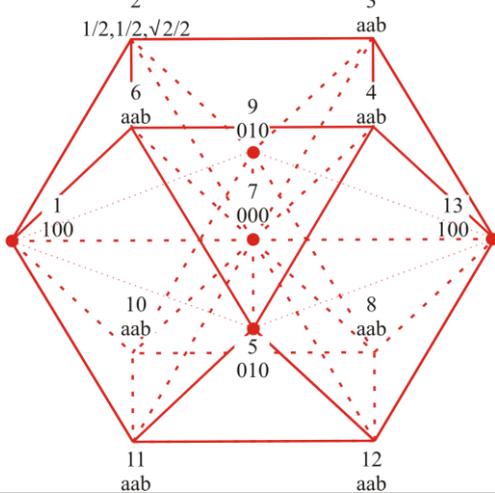 |

Table 2: $2^n$-cubes (left column) and $n$-vector equilibriums (right column).

For $2 \leq n \leq 3$ $n$-vector equilibrium can be oriented in $\mathbb{R}^n$ such that its $2(n-1)$ external vertices and the central vertex define Cartesian coordinate system for $\mathbb{R}^{n-1}$. This cannot be done for $n > 3$ as the number of the remaining external vertices is

$$|v(n)| - 2(n-1) = n^2 - n + 2$$

while $2^n$ is required. $n^2 - n + 2 = 2^n$ only for $1 \leq n \leq 3$.

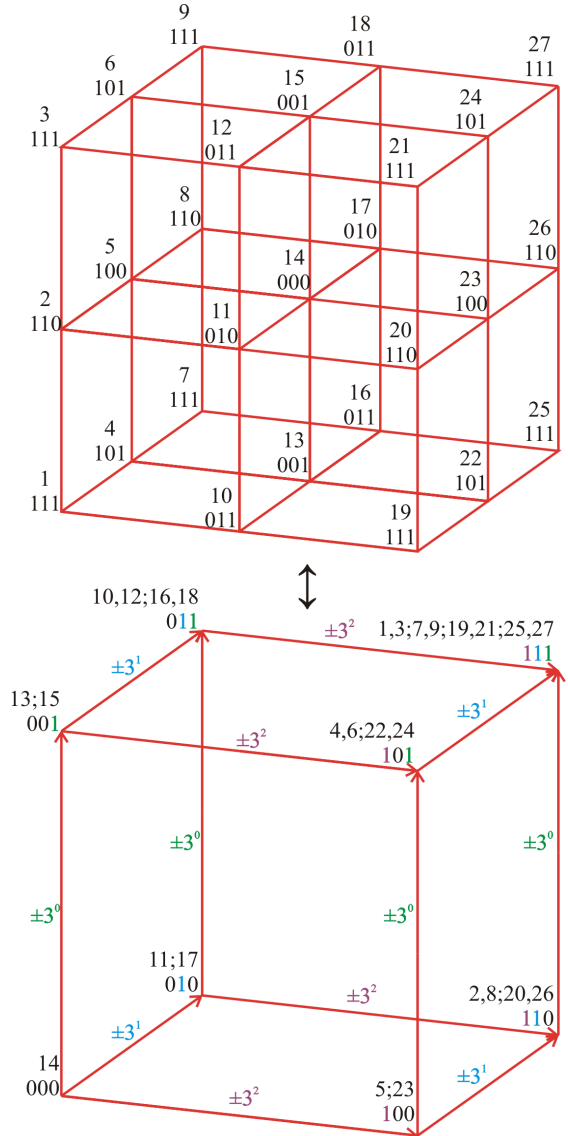

Fig. 9: Relation between $2^n$-cube and $n$-cube ($n = 3$).

$2^n$-cubes follow generalised principles (that is the rules that hold true without exception, according to Buckminster Fuller definition) with regard to their dimensionality, although they do not provide optimal energetic conditions. Unit spheres, for example, placed in vertices of $2^3$-cube will not form a spatially optimal and stable arrangement.

$n$-vector equilibriums, on the other hand, seem to provide optimal energetic conditions. They solve the kissing number problem for $1 \leq n \leq 3$ but they do not follow generalised principles. Synergetic properties of these structures were also recognized by Pope Francis (at least for $n = 3$): "our model must be that of a polyhedron, in which the value of each individual is respected, where «the whole is greater than the part, but it is also greater than the sum of its parts»" and "our model is not the sphere, which is no greater than its parts, where every point is equidistant from the centre, and there are no differences between them. Instead, it is the polyhedron, which reflects the convergence of all its parts, each of which preserves its distinctiveness." [PF20] (pt. 145) and [PF13] (pt. 236).

## V. Discussion

Everyday experience suggests that $n = 4$. Three bidirectional *spatial* dimensions and one unidirectional time.

Biological evolution is a change in the heritable characteristics of biological populations of individuals over successive generations by means of natural selection. These characteristics is the information (genes) passed on from an individual parent(s) to an individual offspring during reproduction. Evolution is the process by which traits that enhance survival and reproduction become more common in successive generations of a population. It is commonly considered as a *self-evident* mechanism because it necessarily follows from three simple facts:
1. (Phenotypic variation) variation of phenotypic traits exists within populations of individuals;
2. (Differential fitness) different traits confer different rates of survival and reproduction; and
3. (Heritability of fitness) these traits can be passed from generation to generation.

Individual is a living entity. A cell, an organism or a neural network in particular. But not a virus or DNA.

Perception[15] is a mapping between an external information (belonging to the individual's *umwelt*) and corresponding memorized information stored in the individual's memory. Sparse distributed memory is ideal for storing a predictive model of the world [PK88].

It has been demonstrated [GR14], however, that the process of memorising information does not require neural networks. Memory has evolved simply to enable the reproductive fitness. Michio Kaku defines the consciousness of an individual simply as the number of feedback loops required to create a model of itself in relationship to space, other individuals, and time [MK15]. It is memory that enables us to perceive *movement* despite Zeno's paradoxes of motion.

Assume that any piece of memorized information is bounded by implicational constraints (ICT) on atomic predicates of an individual's SDM and thus corresponds to a differentiable structure that the individual perceives using Stokes' theorem over its $n-1$ dimensional boundary of perception.

Only for $\mathbb{R}^4$ there exists an uncountable family of non-diffeomorphic differentiable structures which are homeomorphic to $\mathbb{R}^4$ [CT87] (for every such smooth structure there exists a simplicial triangulation, its discrete version [KS77][16]). This property is known as exotic $\mathbb{R}^4$.

---
[15] In a way perception is a synonym for a measurement.

For $n = 1, 2, 3$ any smooth manifold homeomorphic to $\mathbb{R}^n$ is also diffeomorphic to $\mathbb{R}^n$. For $n > 4$ examples of homeomorphic but not diffeomorphic pairs on spheres have been found but are countable (cf. Milnor's sphere for $n = 7$). For $n \neq 4$ exotic $\mathbb{R}^n$'s do not exist [RG83].

This feature of 4-dimensional space, that we experience, indicates that it may be some kind of a prerequisite or $4^{th}$ fact of the biological evolution.

Indeed for $n = 1, 2, 3$ any differentiable structure perceived by an individual would be diffeomorphic to the corresponding differentiable structure that this individual already memorised. There would be only one equivalence class between them. This in turn would have contradicted the principles of biological evolution: no variations of traits would have existed, the same traits would have conferred the same rates of survival and reproduction, and there would have been no need to pass these same traits from generation to generation.

That suggests that only four dimensions allow for variations of traits between any two individuals that perceived the same differentiable structure. If the perception of the organism was statistically tested in the same conditions and in a large number of trials there is no reason to expect that the same equivalence class between perceived and memorised differentiable structure would be applied by each individual and in each trial.

$n > 4$ provides examples of homeomorphic but not diffeomorphic pairs of differentiable structures. But they are countably finite, so a sufficiently large population of individuals would have soon saturated this set and the evolution would have terminated.

Biological evolution would not be possible also if each differentiable structure perceived by an organism could not be diffeomorphic to the corresponding memorised structure. In this case no implicational constraints (ICT) could be formed in an individual's SDM and no meaningful semiosis would be possible between individuals.

If an individual is an observer then homeomorphic but not diffeomorphic pairs of differentiable structures correspond to observer-dependent facts, while homeomorphic and diffeomorphic pairs correspond to observer-independent facts.

Observer-independence has been rejected in a quantum photonic experiment [PEA19] (implementing the gedanken experiment proposed in [CB15]), which demonstrated that no general framework exists in which all observers can reconcile all their recorded facts. That means that there is no (single or unique) objective reality that an observer could communicate to another one.

This by no means boils down to subjectivism since the existence of observer-dependent (aka subjective) facts, as such, does not preclude an existence of such an observer-independent, general framework. But this contradicts the results of this experiment.

Both kinds of facts exist (observer-dependent and observer-independent) and no consistent objective reality can be constructed just by the independent ones.

Observer-dependent facts correspond to non-orthogonal quantum states, while observer-independent to orthogonal ones. The latter can be copied using a single unitary operation.

Also the Ugly duckling theorem implies rejection of observer-independent (aka objective) reality that cannot be constructed solely by observer-independent facts, i.e. equally similar (that is, the same) *objects*, *particles*, etc. This theorem holds trivially for points (not objects) in a space; any two points, insofar as they are distinguishable, are clearly equally similar. For objects, it invalidates the identity of indiscernibles ontological axiom, which is also the $1^{st}$ axiom of metric[17].

The *corollary* of this theorem (assigning individual weights to the predicates (CPP) to assert the similarity of the objects) is just the $2^{nd}$ fact of evolution (differential fitness). Some individuals do it locally better, some do it locally worse, and timeless optimum can never be achieved. We can never *hear the shape of a drum* [MK66]. One simply learns to discern.

For $n = 4$ every individual memorises its own, unique version of observer-dependent reality that it perceives through 3-dimensional boundary of perception. Evolution *creates* the space as we perceive it. And our perception is clearly different to a perception of a Mantis shrimp, a Chameleon or a *Cuscuta* plant for example.

The process of assigning spatial coordinate to a point clearly takes place as we learn to see. Absorption of a photon of visual light by an electron in a cone cell of the eye leads to photoisomerization of retinaldehyde and triggers signal transduction cascade finalised in perception of the information encoded in this photon in a brain. This question was addressed by Časlav Brukner: „one could imagine an experiment on a person whose nerve fibers behind the retina are disconnected and again reconnected at different, randomly chosen, nerve extensions connecting to the brain. It seems reasonable to assume that the neighbouring points of the *object* that is illuminated with light and observed by the person's eye will no more be perceived by the person as neighbouring points. One may wonder if, in the course of further interaction with the environment, the person's brain will start to make sense out of the seen «disordered *classical world*», or if it will post-process the signals to search for more «ordered» structures as a prerequisite for making sense out of them. The latter may eventually nullify the effect of the random reconnection of nerves, and the person will again perceive *the ordinary classical world*." [CB15] (p. 8,9; emphasis added).

The mere fact that we use three-dimensional kinematic models with time to describe the perceived space, move within it, and create numerous inventions that facilitate

---

[16] Thanks for a Moishe Kohan for providing a clarification of this issue on Mathematics Stack Exchange.

[17] Łukaszyk-Karmowski distance function does not follow this axiom.

this movement, does not imply that this space exists, as such.

Taking the above considerations into account, nothing precludes the author to treat the universe (including author's material part) as a topological graph having certain intrinsic properties reflecting the 2$^{nd}$ law of thermodynamics, rather than as something that for *velocities* much lower than the speed of light and for *distances* much higher than Planck length should approximate to *old good classical reality*.

## VI. Acknowledgements

As always, I truly thank my wife for her support, and Mirek, Jacek and Jarek for a lot of inspiring discussions.

| | |
|---|---|
| [AE87] | Amy C. Edmondson, "A Fuller Explanation: The Synergetic Geometry of R. Buckminster Fuller". DOI 10.1007/978-1-4684-7485-5, 1987. |
| [AH03] | Anil N. Hirani, "Discrete Exterior Calculus". PhD Thesis, Caltech 2003. |
| [CB15] | Časlav Brukner, "On the quantum measurement problem". In: Bertlmann R., Zeilinger A. (eds) Quantum [Un]Speakables II. The Frontiers Collection. Springer, Cham. https://doi.org/10.1007/978-3-319-38987-5_5, July 19, 2015. |
| [CB76] | A. Cantoni, P. Butler, "Eigenvalues and Eigenvectors of Symmetric Centrosymmetric Matrices". Linear Algebra and its Applications 13, 275-288, American Elsevier Publishing Company, Inc., 1976. |
| [CL72] | Charles L. Lawson, "Transforming triangulations". Discrete Mathematics, vol 3, Issue 4, 1972. |
| [CT87] | Clifford Henry Taubes, "Gauge theory on asymptotically periodic {4}-manifolds". J. Differential Geom., Volume 25, Number 3 (1987), 363-430. |
| [DKT08] | Mathieu Desbrun, Eva Kanso, Yiying Tong, "Discrete Differential Forms for Computational Modeling". Discrete Differential Geometry. Oberwolfach Seminars, vol 38. Birkhäuser Basel (2008). |
| [GR14] | Monica Gagliano, Michael Renton, Martial Depczynski, Stefano Mancuso, "Experience teaches plants to learn faster and forget slower in environments where it matters". Oecologia 175.1 (2014): 63-72. |
| [KC19a] | Keenan Crane, "The *n*-dimensional cotangent formula". Online note. URL: https://www.cs.cmu.edu/~kmcrane/Projects/Other/nDCotanFormula.pdf, 2019. |
| [KS77] | Robion C. Kirby, Laurence C. Siebenmann, "Foundational Essays on Topological Manifolds, Smoothings, and Triangulations". Princeton University Press and University of Tokyo Press, Princeton, New Jersey, 1977. |
| [LB77] | Ludwig Boltzmann, "Über die Beziehung zwischen dem zweiten Hauptsatze des mechanischen Wärmetheorie und der Wahrscheinlichkeitsrechnung, respective den Sätzen über das Wärmegleichgewicht" Von dem c. M. Ludwig Boltzmann in Graz Sitzb. d. Kaiserlichen Akademie der Wissenschaften, mathematich-naturwissen Cl. LXXVI, Abt II, 1877, pp. 373-435., av. at http://users.polytech.unice.fr/~leroux/boltztrad.pdf |
| [MK15] | Michio Kaku, "The Future of the Mind". ISBN-13: 978-0385530828, 25 Feb 2014. |
| [MK66] | Mark Kac. "Can one hear the shape of a drum?" The American Mathematical Monthly, 73(4):1–23, 1966. |
| [MW17] | Max Wardetzky, "Generalized Barycentric Coordinates in Computer Graphics and Computational Mechanics". Chapter 5, A Primer on Laplacians, CRC Press, 2017. |
| [PDEA20] | Peter B. Denton, Stephen J. Parke, Terence Tao, and Xining Zhang, "Eigenvectors from Eigenvalues: a Survey of a Basic Identity in Linear Algebra". TBD, https://www.osti.gov/biblio/1561548, 4 Mar 2020. |
| [PEA19] | Massimiliano Proietti, Alexander Pickston, Francesco Graffitti, Peter Barrow, Dmytro Kundys, Cyril Branciard, Martin Ringbauer, Alessandro Fedrizzi, "Experimental test of local observer independence". Science Advances, Vol. 5, no. 9, eaaw9832, DOI: 10.1126/sciadv.aaw9832, 20 Sep 2019. |
| [PF13] | Pope Francis, Apostolic Exhortation "Evangelii Gaudium". 24 November 2013. |
| [PF20] | Pope Francis, "Fratelli tutti". 3 Oct 2020. |
| [PK88] | Pentti Kanerva, "Sparse Distributed Memory". The MIT Press. ISBN 978-0-262-11132-4 1988. |
| [RCEA09] | Renjie Chen, Yin Xua, Craig Gotsman, Ligang Liu, "A spectral characterization of the Delaunay triangulation". Computer Aided Geometric Design, Elsevier Vol. 27, Issue 4, May 2010. |
| [RG83] | Robert E. Gompf, "Three Exotic R$^4$'s and other Anomalies". J. Differential Geometry 18 (1983) 317-328. |
| [SL03] | Szymon Łukaszyk, "A new concept of probability metric and its applications in approximation of scattered data sets". Computational Mechanics, 33-4, 299–304, doi:10.1007/s00466-003-0532-2 (2003). |
| [SW69] | Satosi Watanabe, "Knowing and Guessing: A Quantitative Study of Inference and Information". New York: Wiley. ISBN 0-471-92130-0. LCCN 68-56165. (1969). |
| [SW86] | Satosi Watanabe, "Epistemological Relativity - Logico-Linguistic Source of Relativity". Annals of the Japan Association for Philosophy of Science, March 1986. |
| [TK78] | Thomas S. Kuhn, "Black-Body Theory and the Quantum Discontinuity, 1894-1912". Oxford University Press, 1978. |
| [VNB08] | Péter Várlaki, László Nádai, József Bokor, "Number Archetypes and 'Background' Control Theory Concerning the Fine Structure Constant". Acta Polytechnica Hungarica Vol. 5, No. 2, 2008. |

MATLAB scripts that were created to generate matrices of the graphs discussed in this paper, as well as an Excel spreadsheet containing their particular forms and the sequences of their eigenvalues are available at https://github.com/szluk/FourCubes.